\documentclass[11pt,leqno]{amsart}
\usepackage{amsmath,amssymb,amsthm}
\usepackage{hyperref}

\usepackage{bbm}

\newcommand{\re}{\operatorname{Re}}

\renewcommand{\geq}{\geqslant}
\renewcommand{\leq}{\leqslant}

\newcommand{\1}[1]{\operatorname{\textbf{1}}}

\newtheorem{theorem}{Theorem}[section]

\theoremstyle{definition}

\theoremstyle{remark}
\newtheorem{remark}[theorem]{Remark}
\numberwithin{equation}{section}

\def\fnote#1{\footnote}

\def\ignora#1{}
\def\n3#1{\left\vert  \! \left\vert \! \left\vert \, #1 \, \right\vert \!
  \right\vert \! \right\vert }

\newcommand{\iten}{\ensuremath{\widehat{\otimes}_\varepsilon}}
\newcommand{\pten}{\ensuremath{\widehat{\otimes}_\pi}}

\begin{document}

\title{ Stability results of octahedrality in tensor product spaces }

\author{ Abraham Rueda Zoca }\address{Universidad de Granada, Facultad de Ciencias. Departamento de An\'{a}lisis Matem\'{a}tico, 18071-Granada
(Spain)} \email{ abrahamrueda@ugr.es}
\urladdr{\url{https://arzenglish.wordpress.com}}



\maketitle

\markboth{ABRAHAM RUEDA ZOCA}{STABILITY RESULTS OF OCTAHEDRALITY IN TENSOR PRODUCT SPACES}

\begin{abstract}
We prove that there exists a finite-dimensional Banach space $X$ such that $L_1^\mathbb C([0,1])\iten X$ fails the strong diameter two property and $L_\infty^\mathbb C([0,1])\pten X^*$ fails to have octahedral norm. This proves that the octahedrality of the norm (respectively the strong diameter two property) is not automatically inherited from one factor by taking projective tensor product (respectively injective tensor product), which answers \cite[Question 4.4]{llr}.
\end{abstract}

\section{Introduction}

A Banach space $X$ is said to have an \textit{octahedral norm} if, for every finite-dimensional subspace $Y$ of $X$ and every $\varepsilon>0$, there exists $x\in S_X$ such that the inequality
$$\Vert y+\lambda x\Vert\geq (1-\varepsilon)(\Vert y\Vert+\vert\lambda\vert)$$
holds for every $y\in Y$ and every $\lambda\in\mathbb R$. Octahedral norms were intensively studied at the end of the eighties \cite{dev, godefroy,gk}. Recently, continuing with a research line started by R. Deville (see Proposition 3. and Remark (c) in \cite{dev}), the authors of \cite{blr14} gave a characterisation of octahedral norms in terms of the presence of diameter two properties in the dual space (see \cite{blr14} for details). See also \cite{hlp2,lan} for further nice results and different reformulations of octahedral norms. The characterisation given in \cite{blr14} has motivated a lot of new research on octahedral norms and, in particular, a lot of new result about octahedrality in tensor product spaces have appeared in the last few years (see \cite{blr,blr3,hlp,lan,lr2020,llr,mrpre} and references therein).

An open problem posed in this line \cite[Question 4.4]{llr} is the following: given two Banach spaces $X$ and $Y$, does $X\pten Y$ have an octahedral norm whenever $X$ has an octahedral norm? There are some examples where the answer is affirmative such as when $X=L_1(\mu)$ (or more generally, when $X$ is a non-reflexive $L$-embedded space such that $X^{**}$ has the metric approximation property \cite[Theorem 4.3]{llr2}), when the Cunningham algebra of $X$ is infinite-dimensional, when $X$ is the dual of an ASQ space (say $X_*$) and $Y$ is a dual space (say of $Y_*$) such that $X\pten Y=(X_*\iten Y_*)^*$ (see \cite[Section 4]{llr}) or when $X=\mathcal F(M)$ for any infinite metric space $M$ which is not uniformly bounded \cite[Theorem 4.1]{mrpre}. However, to the best of our knowledge, the above mentioned question \cite[Question 4.4]{llr} still remains open. The aim of this short note is to give a negative answer to \cite[Question 4.4]{llr} by proving that, in general, the projective tensor product does not preserve octahedrality from one of its factors. The counterexample exposed in Theorem~\ref{theo:contraproj} is strongly motivated by an example presented by V. Kadets, N. J. Kalton and D. Werner \cite[Corollary 4.3]{kkw}, where the authors exhibited an example of a projective tensor product failing the Daugavet property even though one of its factors of such space does have the Daugavet property.

\section{Notation and preliminary results}\label{section:notation}

For the notation section we consider real or complex Banach spaces. Given a Banach space $X$, we denote by $B_X$ and $S_X$ the closed unit ball and the closed unit sphere of $X$. We also denote by $X^*$ the topological dual of $X$.

By a slice of the unit ball we denote a set of the form
$$S(B_X,x^*,\alpha):=\{x\in B_X: \re x^*(x)>\sup \re(x^*(B_X))-\alpha\},$$
where $x^*\in X^*$ and $\alpha>0$. We also consider a \textit{convex combination of slices of $B_X$} a set of the following form
$$\sum_{i=1}^n\lambda_i S_i,$$
where $S_1,\ldots, S_n$ are slices of $B_X$ and  $\lambda_1,\ldots, \lambda_n\in [0,1]$ satisfy that $\sum_{i=1}^n\lambda_i=1$.

We say that $X$ has the \textit{strong diameter two property (SD2P)} if every convex combination of $B_X$ has diameter exactly two. It is known that a Banach space $X$ has the SD2P if, and only if, the norm of $X^*$ is octahedral \cite[Corollary 2.2]{blr14}.

Another characterisation of the SD2P is the following: a Banach space $X$ has the SD2P if, and only if, given any convex combination of slices $C$ of $B_X$ and any $\varepsilon>0$ there exists $x\in C$ with $\Vert x\Vert>1-\varepsilon$ \cite[Theorem 3.1]{lmr}.

Given a complex Banach space $X$ we say that $X$ is \textit{complex uniformly convex} if there exists a function $\delta:\mathbb R^+\longrightarrow \mathbb R^+$ with $\delta(\varepsilon)\rightarrow 0\ (\varepsilon\rightarrow 0)$ and with the property that
$$\left.
\begin{array}{c}
1-\varepsilon<x\leq 1, y\in B_X  \\
\Vert x\pm y\Vert\leq 1,\\
\Vert x\pm i y\Vert\leq 1
\end{array}
\right\}\Rightarrow \Vert y\Vert<\delta(\varepsilon).$$
Observe that the space $L_1^\mathbb C([0,1])$, the space of all complex-valued $L_1$ functions, is complex uniformly convex (indeed, $\delta(\varepsilon)=\sqrt{\varepsilon}(4+2\sqrt{1+2\varepsilon})$ can be taken, see  \cite[Theorem 1]{glov}).

Given two Banach spaces $X$ and $Y$ we will denote by $L(X,Y)$ (resp. $K(X,Y)$) the space of all linear and bounded (resp. compact linear) operators from $X$ to $Y$, and we will denote by $X\pten Y$ and $X\iten Y$ the projective and the injective tensor product of $X$ and $Y$. See \cite{rya} for a detailed treatment of the tensor product theory.

\section{Main result}

Now we present the result of the paper.

\begin{theorem}\label{theo:sd2pnoinje}
There exists a two dimensional complex Banach space $X$ such that $L_1^\mathbb C([0,1])\iten X$ fails the strong diameter two property, where $L_1^\mathbb C([0,1])$ stands for the space of complex-valued $L_1$-functions.
\end{theorem}

The idea of the proof below is strongly inspired by the proof of \cite[Theorem 4.2]{kkw}.

\begin{proof} Write for simplicity $L_1=L_1^\mathbb C([0,1])$. Let $X$ the subspace of complex $\ell_\infty^{10}$ spanned by the vectors $x_1:=(1,1,1,1,1,\frac{1}{2},-\frac{1}{2}, \frac{i}{2},-\frac{i}{2},0)$ and\\  $x_2:=(0,\frac{1}{2},-\frac{1}{2}, \frac{i}{2},-\frac{i}{2},1,1,1,1,1)$. Then $L_1\iten X$ identifies with the space of $10$-tuples $(f_1,f_2,\ldots, f_{10})$, with $f_i\in L_1$, of the form $g_1\otimes x_1+g_2\otimes x_2$ and with the norm $\Vert f\Vert:=\max\limits_{1\leq k\leq 10}\Vert f_k\Vert$ (this identification comes from the identification $L_1\iten \ell_\infty^{10}=\ell_\infty^{10}(L_1)$ \cite[Chapter 3.2]{rya} and by the well known fact that the injective tensor product respects subspaces isometrically \cite[Proposition 3.2]{rya}).

Let us prove that $L_1\iten X$ fails the SD2P. In order to do so, given $\varepsilon>0$, consider
$$S_\varepsilon:=\left\{f=(f_1,\ldots, f_{10})\in L_1\iten X: \Vert f\Vert\leq 1, \re \int_0^1 f_1>1-\varepsilon\right\},$$
which is a slice of $B_X$. Given any $g_1\otimes x_1+g_2\otimes x_2\in S_\varepsilon$ it follows
\[\begin{split}
1\geq \Vert g_1\Vert\geq \re \int_0^1 g_1(t)\ dt>1-\varepsilon\\
\max\left\{ \left\Vert g_1\pm \frac{g_2}{2}\right\Vert,\left\Vert g_1\pm i \frac{g_2}{2}\right\Vert\right\}\leq 1,
\end{split}
\]
where the above conditions hold using that $\Vert g_1\otimes x_1+g_2\otimes x_2\Vert\leq 1$ and having a look to the five first coordinates. Since $L_1$ is complex uniformly convex we infer $\Vert g_2\Vert\leq \delta(\varepsilon)$, where $\delta(\varepsilon)$ is a positive function on $\varepsilon$ with $\lim\limits_{\varepsilon\rightarrow 0} \delta(\varepsilon)=0$ (see Section~\ref{section:notation} and the references given there). 

On the other hand, consider the slice
$$T_\varepsilon:=\left\{f=(f_1,\ldots, f_{10})\in L_1\iten X: \Vert f\Vert\leq 1, \re \int_0^1 f_{10}>1-\varepsilon\right\}$$
By a similar argument of complex uniform convexity, given $h_1\otimes x_1+h_2\otimes x_2\in T_\varepsilon$ we derive that $\Vert h_1\Vert\leq \delta(\varepsilon)$.

Now, given $z_1=g_1\otimes x_1+g_2\otimes x_2\in S_\varepsilon$ and $z_2=h_1\otimes x_1+h_2\otimes x_2\in T_\varepsilon$ we get
$$z_1+z_2=(g_1+h_1)\otimes x_1+(g_2+h_2)\otimes x_2,$$
where $\max\{\Vert g_2\Vert, \Vert h_1\Vert\}\leq \delta(\varepsilon)$. Hence
\[\begin{split}
\Vert z_1+z_2\Vert& =\max\left\{\Vert g_1+h_1\Vert, \left\Vert (g_1+h_1)\pm \frac{g_2+h_2}{2}\right\Vert\right. ,\\
& \left\Vert (g_1+h_1)\pm i\frac{g_2+h_2}{2}\right\Vert, \left\Vert (g_2+h_2)\pm \frac{g_1+h_1}{2}\right\Vert,\\
& \left. \left\Vert (g_2+h_2)\pm i\frac{g_1+h_1}{2}\right\Vert,\left\Vert g_2\right\Vert\right\}\\
& \leq \max\left\{ \Vert g_1\Vert+\Vert h_1\Vert+\frac{\Vert g_2\Vert+\Vert h_2\Vert}{2}, \Vert g_2\Vert+\Vert h_2\Vert+\frac{\Vert g_1\Vert+\Vert h_1\Vert}{2}\right\}\\
& \leq 1+\delta(\varepsilon)+\frac{1+\delta(\varepsilon)}{2}=\frac{3}{2}+2\delta(\varepsilon).
\end{split}\]
The arbitrariness of $z_1$ and $z_2$ implies that $S_\varepsilon+T_\varepsilon\subseteq \left(\frac{3}{2}+2\delta(\varepsilon)\right)B_X$.

Now consider $C_\varepsilon:=\frac{S_\varepsilon+T_\varepsilon}{2}$, which is a convex combination of slices. The above estimate implies that $C_\varepsilon\subseteq (\frac{3}{4}+\delta(\varepsilon))B_X$. If we select $\varepsilon$ small enough to get $\delta(\varepsilon)<\frac{1}{4}$, we get $C_\varepsilon$ contained in ball centred at $0$ of radius strictly smaller than $1$. Now \cite[Theorem 3.1]{lmr} implies that $L_1\iten X$ fails the SD2P.
\end{proof}

Since it is known that a Banach space $X$ has the SD2P if, and only if, the norm of $X^*$ is octahedral \cite[Corollary 2.2]{blr14} we have the following result.

\begin{theorem}\label{theo:contraproj}
There exists a two dimensional complex Banach space $Y$ such that the norm of $L_\infty^\mathbb C([0,1])\pten Y$ is not octahedral, where $L_\infty^\mathbb C([0,1])$ stands for the space of complex-valued $L_\infty$-functions.
\end{theorem}

 \begin{proof}
Take the Banach space $X$ from Theorem~\ref{theo:sd2pnoinje} such that $L_1\iten X$ fails the SD2P. Taking $Y=X^*$ we get that the norm of $(L_1\iten X)^*=L_\infty\pten X^*=L_\infty\pten Y$ is not octahedral.
 \end{proof}

Let us now observe some consequences of the example exposed here:

\begin{remark}
\begin{enumerate}
\item The previous result gives a negative answer to \cite[Question 4.4]{llr}, where it is wondered whether octahedrality is preserved by just one factor by taking projective tensor product. Indeed, the space $L_\infty^\mathbb C([0,1])$ has the Daugavet property (see e.g. \cite[Section 2, Examples (b)]{werner}), which implies that its norm is octahedral (see \cite[Lemma 2.8]{kssw}), but Theorem~\ref{theo:contraproj} implies that $L_\infty^\mathbb C([0,1])\pten Y$ even fails to have octahedral norm.

\item Given two Banach spaces $X$ and $Y$ there are some stability results of the SD2P in $L(X,Y)$ assuming hypothesis about SD2P just in $X$ (see. e.g. \cite[Corollary 3.7]{abl}). However, Theorem~\ref{theo:sd2pnoinje} implies that the SD2P is not preserved in general from one of the factors by taking space of linear and continuous operators. Indeed, taking the space $X$ of Theorem~\ref{theo:sd2pnoinje}, the space $L_1^\mathbb C([0,1])\iten X=K(L_\infty^\mathbb C([0,1]),X)$ fails the SD2P, consequently its bidual $(L_1^\mathbb C([0,1])\iten X)^{**}=(L_\infty^\mathbb C([0,1])\pten X^*)^*=L(L_\infty([0,1]),X^{**})$ fails the SD2P too (c.f. e.g. \cite[Proposition 2.14]{lan}). However, $L_\infty^\mathbb C([0,1])$ has the SD2P because $L_\infty^\mathbb C([0,1])$ has the Daugavet property (see e.g. \cite[Section 2, Examples (b)]{werner}) and Daugavet spaces have SD2P (c.f. e.g. \cite[Theorem 4.4]{aln}).
\end{enumerate}
\end{remark}

\section*{Acknowledgements}  

This work was supported by MCIN/AEI/10.13039/501100011033: grant PID2021-122126NB-C31, Junta de Andaluc\'ia: grant FQM-0185, by Fundaci\'on S\'eneca: ACyT Regi\'on de Murcia: grant 21955/PI/22 and by Generalitat Valenciana: grant CIGE/2022/97.

\end{document}